\documentclass{article}
\usepackage{graphicx} 
\usepackage{parskip}
\usepackage{amsmath}
\usepackage[colorlinks=true, linkcolor=blue, urlcolor=blue, citecolor=blue]{hyperref}

\title{A Summation-Based Algorithm For Integer Factorization}
\author{Justin Friedlander}
\date{March 2025}

\begin{document}

\maketitle

\section{Introduction}
Numerous methods have been considered to create a fast integer factorization algorithm. Despite its apparent simplicity, the difficulty to find such an algorithm plays a crucial role in modern cryptography, notably, in the security of RSA encryption. Some approaches to factoring integers quickly include the Trial Division method, Pollard's Rho and p-1 methods, and various Sieve algorithms \cite{wagstaff2002}.
\\\\
This paper introduces a new method that converts an integer into a sum in base-2. By combining a base-10 and base-2 representation of the integer, an algorithm on the order of $\sqrt{n}$ time complexity can convert that sum to a product of two integers, thus factoring the original number.

\section{Method}
\subsection*{Step One: Iterating Through $j$ and $i$}
Let $n=pq$ for integers $n$, $p$, and $q$. Note that $p$ and $q$ can be written in base-2. Consider, however, the highest power of $p$ and $q$. That is, $\lfloor log_2(p) \rfloor$ and $\lfloor log_2(q) \rfloor$. WLOG, let $p \geq q$. Let $j = \lfloor log_2(p) \rfloor$ and $i = \lfloor log_2(q) \rfloor$. Note that $p = 2^j + c_i$ and $q = 2^i + c_j$ for some integers $c_i < 2^j$ and $c_j < 2^i$.
\\\\
Note that now $n=pq=(2^j+c_i)(2^i+c_j)=2^{j+i}+c_j2^j+c_i2^i+c_jc_i$.
\\\\
We can also represent $n$ in base-2, however, it may or may not be identical to our representation of $pq$.

\begin{quote}
    \itshape Theorem 1:
    \begin{quote}
        Let $n=2^k+c_k$ for $k = \lfloor log_2(n) \rfloor$ and $c_k < 2^k$.
        \\
        Claim: $k = j + i$ or $k = j + i + 1$ for all $n$, $p$, and $q$.
        \\\\
        Proof:\\
        Lower Bound - \ $n = 2^k + c_k = (2^j + c_i)(2^i + c_j)$. Let $c_j = c_i = 0$. Now, $n = 2^k + c_k = 2^{j+i}.$ Since $k$ is the largest power of $2$ before increasing above $n$, $j+i = k$. Thus $j + i \leq k$ for any arbitrary $c_j$ and $c_i$.
        \\\\
        Upper Bound - \ $c_j < 2^i$ and $c_i < 2^j$. Thus, 
        \begin{align*}
            n &= 2^k + c_k \\
              &= (2^j+c_i)(2^i+c_j)\\
              &= 2^{j+i}+c_j2^j+c_i2^i+c_jc_i\\
              &< 2^{j+i} + 2^i2^j + 2^j2^i + 2^i2^j\\
              &= 4*2^{j+i}
        \end{align*}
        Since $2^k + c_k < 4*2^{j+i}$, then $2^k + c_k < 2^{j+i+2}$. Thus, to get the left-hand-side and right-hand-side to be equal, we must decrement the right hand side by at least one. This leaves $2^k + c_k = 2^{j+i+1} + c_{decrement}$. Again, since $2^k$ is the largest power of $2$ before increasing above $n$, $k = j + i + 1$.\\\\
        
        Thus, $k \leq j + i + 1$, so $j + i \leq k \leq j + i + 1$ for all $n$, $p$, and $q$.
        
    \end{quote}
\end{quote}

The implications of Theorem 1 are that the algorithm will have to run once to check the case where $k = j + i$, and a second time to check if $k = j + i + 1$ in the worst case scenario.
\\\\
Additionally, when given a power $k$, the numbers $j$ and $i$ are unknown. Thus, the algorihtm must search through all combinations of $j$ and $i$ such that $k = j + i$ or $k = j + i + 1 $.

\subsection*{Step Two: Iterating Through $c_J$}
Since we are iterating over all combinations of $j$ and $i$, for this next part of the the algorithm, we can assume our choices of $j$ and $i$ are the correct choices that correspond with $p$ and $q$. That is, $j = \lfloor log_2(p) \rfloor$ and $i = \lfloor log_2(q) \rfloor$. Since the following argument is nearly identical for $k = j + i$ and $k = j + i + 1$, we will assume $k = j + i$ for simplicity.
\\\\
We know $n = 2^k + c_k = 2^{j+i}+c_j2^j+c_i2^i+c_jc_i$ and $2^k = 2^{j+i}$. Thus, $c_k = c_j2^j + c_i2^i + c_jc_i$. We can represent $c_k$ in this form by reducing it in base-2. Here is an example of such a process:
\begin{quote}
    $c_k = 61$, $j=4$, $i=2$.\\
    $c_k = 2^5 + 2^4 + 2^3 + 2^2 + 2^0 = 2*2^4 + 2^4 + 2*2^2 + 2^2 + 2^0 = 3*2^4 + 3*2^2 + 1 = 3*2^j + 3*2^i + 1$
\end{quote}
We can define $c_J$ and $c_I$ to equal the respective coefficients of $2^j$ and $2^i$, and $B$ to equal the coefficient of $2^0$ after reducing $c_k$ to this form. Notice that $c_J2^j + c_I2^i + B = (c_J-e)2^j + (c_I + e2^{j-i})2^i + B = c_J2^j + (c_I-d)2^i + (B + d2^i)$ for some integers $e$ and $d$. From the above example, we can write:
\begin{quote}
    $3*2^4 + 3*2^2 + 1 = (3-2)*2^4 + (3 + 2*2^{4-2})2^2 + 1 = 2^4 + (11 - 2)2^2 + (1 + 2*2^2) = 2^4 + 9*2^2 + 9$
\end{quote}
Now, if we re-introduce the $2^k$ term, we get
\begin{quote}
    $2^k + c_k = 2^{4+2} + 2^4 + 9*2^2 + 9 = (2^4 + 9)(2^2 + 1) = 25 * 5 = pq$
\end{quote}
Notice that we will know we have achieved the correct coefficients for $c_j$ and $c_i$ when $c_jc_i = b$ where $b$ is our $2^0$ coefficient.
\\\\
The algorithm I have found that converts from $c_J$, $c_I$, and $B$ to $c_j$, $c_i$, and $b$ must consider, in the worst case, all the iterations of the $2^j$ coefficient from $c_J$ to $1$. Since we are iterating through all values of this coefficient, we can assume that this coefficient is $c_j$.

\subsection*{Step Three: Finding $c_i$}
Let $e = c_J - c_j$ and $c_I^{\prime} = c_I + e*2^{j-i}$. We can use the equation below to find the difference $d$ between $c_I^{\prime}$ and $c_i$:
\begin{quote}
    Equation 1: $\frac{(c_J - e)c_I^{\prime} + B}{c_J - e + 2^i} = d$
\end{quote}
From this, we can compute $c_i$ from $c_I^{\prime} - d$ and $c_j$ from $c_J - e$. Since we know our $c_j$ and $c_i$, and we know $j$ and $i$, we know the term $(2^j + c_i)(2^i + c_j) = n$, so we can deduce our $p$ and $q$.

\section{Time Complexity}
In the first part of the algorithm, we are iterating through all the combinations of $j$ and $i$ such that $k = j + i$ or $k = j + i + 1$. Since $k$ is approximately $log(n)$, this step requires approximately $log(n)$ iterations. In the second step of the algorithm, we must iterate through all the coefficients of the $2^j$ term. Since $c_j < 2^i$, and $j \geq i$, in the worst case we have $c_J \leq \sqrt{n}$. This means that this step can take $\sqrt{n}$ iterations. In the third step, we compute $c_i$ from $c_j$, which is a constant time computation.
\\\\
Thus, the algorithm as a whole seems to take O($\sqrt{n}log(n)$) time to run. Closer inspection, however, reveals one minor improvement to this number. When $j \approx k$, then $i \approx 0$ because $k - j = i$. In this case, $c_J$ is much closer to 0 than $\sqrt{n}$. More generally, each iteration of $j$ and $i$ increases the possible values of $c_J$ by approximately a factor of 2. Thus, the total number of operations performed in this algorithm is closer to $2*\sum_{k=0}^{log_2(n)} \frac{\sqrt{n}}{2^k} \approx 4\sqrt{n}$. Thus, the run-time of this algorithm is on the order of $\sqrt{n}$.

\section{Discussion}
This algorithm falls short of improving upon the time complexity of the General Number Field Sieve \cite{lenstra1993}, but it does introduce a new method to factoring integers that, as far as I am aware, has not previously been considered. After an analysis beyond the scope of this paper, I do not believe it is possible to significantly reduce the time complexity of this algorithm without changing to a new algorithm entirely. Thus, I am now exploring quantum computing options that may open the door to further optimizations, and I encourage others interested in this approach to do the same.
\\
A Python implementation of the classical algorithm can be found \href{https://py3.codeskulptor.org/#user310_4yZpp8XWfN_5.py}{here}.

\end{document}